\documentclass[oneside,openany]{amsart}
\usepackage{amssymb}
\usepackage[frame,ps,dvips,matrix,arrow,curve]{xy}
\input{xypic}
\newcommand{\colim}{\ensuremath{\mathrm{colim}}}
\newcommand{\ilim}{\ensuremath{\mathrm{lim}}}
\newcommand{\hocolim}{\ensuremath{\mathrm{hocolim}}}

\newtheorem{lemma}[equation]{Lemma}
\newtheorem{prop}[equation]{Proposition}
\newtheorem{theor}[equation]{Theorem}

\newtheorem{coro}[equation]{Corollary}

\numberwithin{equation}{section}

\begin{document}
\title{Higher string topology on general spaces}
\author{Po Hu}

\maketitle
 
\section{Introduction}

In~\cite{cs}, Chas and Sullivan considered the free loop space 
$LM = Maps(S^1, M)$ for a smooth orientable compact manifold 
$M$.  They used geometric methods to show that $H_{\ast}LM$ 
has, among other things, the 
structure of a Batalin-Vilkovisky algebra. More precisely, $C_{\ast}LM$
has an
action by the operad of chain complexes $C_{\ast}\mathcal{D}_2$, 
where $\mathcal{D}_2$ is the framed $2$-dimensional 
little disk operad. Each element of 
$\mathcal{D}_2(n)$ is an ordered configuration of $n$ disjoint $2$-dimensional
little disks $d_1, \ldots d_n$ inside the $2$-dimensional disk $D^2$, together with 
a linear parametrization $\lambda_i: S^1 \rightarrow c_i$, where $c_i$
is the boundary of $d_i$, for each $i$. 
Also, for such 
an $M$, there is a chain homotopy equivalence
\begin{equation}
C_{\ast}(LM) \simeq HC^{\ast}(C^{\ast}M)
\end{equation}
the Hochschild cohomology complex of the cochain complex of $M$.

In~\cite{cj}, Cohen and Jones gave a homotopy theoretical realization 
of the result of Chas and Sullivan. Namely, let $\nu_M$ be the 
stable normal bundle of the manifold $M$. One can pull back $\nu_M$ 
along the evaluation map $LM \rightarrow M$, and consider 
the stable bundle $\nu_M$ over $LM$. Then the Thom complex
\[ LM^{\nu_M} \]
is an algebra over the (unframed) $2$-dimensional little disk operad in the 
category of spectra (using the modern theory of strictly commutative,
associative and unital ring spectra).
\vspace{4mm}

In this note, we give a generalized analogue of the results of Chas and 
Sullivan, and of Cohen and Jones.  One direction of the generalization 
is that we consider higher dimensional unframed little disk operads $\mathcal{C}_k$
for all $k$. Also, rather than considering smooth 
orientable manifolds, we consider general finite-dimensional simplicial 
complexes.  
In Section~\ref{topsec} below, we shall construct for each finite 
CW-complex $X$ a spectrum $Maps(S^k, X)^{S(X)}$, and prove the following 
result:

\begin{theor}
For $k \geq 1$ and $X$ as above, 
\[ C_{\ast}Maps(S^k, X)^{S(X)} \]
is naturally homotopy equivalent to
a $\mathcal{C}_{k+1}$-algebra in the category of chain complexes, where 
$C_{\ast}$ denotes chain complexes with coefficients in 
a field $K$.
\label{maintop}
\end{theor}

The method is by showing that $C_{\ast}Maps(S^k, X)^{S(X)}$ is 
naturally equivalent to a Hochschild cohomology complex in 
the category of $\mathcal{C}_k$-algebras, and
use the result of~\cite{hkv} that Hochschild cohomology complexes
in the category of $\mathcal{C}_k$-algebras are $\mathcal{C}_{k+1}$-algebras. 
There is a subtle reason why the present paper is not a direct 
generalization of the Batalin-Vilkovisky structure of~\cite{cs}: although
it is natural to conjecture that the Batalin-Vilkovisky structure 
of~\cite{cs} is the same as the one obtained from Hochschild 
cohomology, as far as I know, at present this has not been proved. 
It should be mentioned that Kallel and Salvatore~\cite{ks} and 
Chataur~\cite{chat} use different, more geometrical, approaches to 
higher-dimensional string topology. Results in this direction 
were also announced by Sasha Voronov and Dennis Sullivan.

In Section~\ref{algsec}, we shall discuss some algebraic aspects
of this story. Namely, we shall introduce the concept of 
Koszul duality on general categories outlined in~\cite{hkv}, 
and related this to the Ginzburg-Kapranov Koszul duality 
for operad~\cite{gk}. Using this technique, we shall prove 
Kontsevich's conjecture that Quillen cohomology in
based $\mathcal{C}_k$-algebras is a shift of their Hochschild
cohomology (Theorem~\ref{shift}). We shall also prove that
the operad $\mathcal{C}_k$ is derived-Koszul 
dual to itself up to a shift (Theorem~\ref{selfdual1}). Therefore, one 
has a natural notion of homotopy Koszul dual $\mathcal{C}_k$-algebras.

Recall that for any $X$, 
the $k$-fold loop space $\Omega^k X$ has the structure of a
$\mathcal{C}_k$-algebra (see~\cite{may}). Namely, suppose that $\alpha \in 
\mathcal{C}_{k}(n)$, so $\alpha$ is a configuration of $n$ $k$-dimensional 
little disks inside $D^k$. Then for $f_1, f_2, \ldots, f_n: S^k \rightarrow X$, 
\[ \alpha(f_1, \ldots, f_n) : S^k \rightarrow X \]
is given by sending the $i$-th little disk in the 
configuration $\alpha$ to $X$ via $f_i$, and sending the rest of 
$D^k$ to the basepoint of $X$. Hence, $C_{\ast}\Omega^k X$ has the structure 
of a $C_{\ast}\mathcal{C}_k$-algebra.
Now one has the following

\begin{prop}
For a space $X$, 
$C^{\ast}X$ and $C_{\ast}\Omega^k X$ are Koszul dual to each other in the 
homotopy 
category of $\mathcal{C}_k$-algebras.
\label{cochain}
\end{prop}

For $k = 1$, this is a part of a conjecture made by R. Cohen 
at the 2002 Aarhus conference.

Finally, in Section~\ref{hochsec}, we shall prove that 
Koszul-dual $\mathcal{C}_k$-algebras have the same (unbased) 
Hochschild cohomologies (Theorem~\ref{kossame}). One should note 
that both based and unbased Hochschild cohomologies of 
$\mathcal{C}_k$-algebras are non-trivial notions rigorously 
introduced in~\cite{hkv}. Returning to spaces, this gives the following

\begin{coro}
For a finite-dimensional $k+1$-connected simplicial complex $X$,
\begin{equation*}
\begin{split}
C_{\ast}Maps(S^k, X)^{S(X)} & \simeq HC^{\ast}_{C_{\ast}\mathcal{C}_k}
(C^{\ast}X) \\
& \simeq HC^{\ast}_{C_{\ast}\mathcal{C}_k}(C_{\ast}\Omega^k X).
\end{split}
\end{equation*}
\label{maincoro}
\end{coro}

\begin{proof}
The first equivalence is from Theorem~\ref{firsteq} below.
The second equivalence follows from 
Theorem~\ref{kossame} and Proposition~\ref{cochain}. 
\end{proof}

For $k =1$, this was proved by Y. F\'{e}lix, L. Menichi and J. Thomas~\cite{fmt}.

Note that Corollary~\ref{maincoro} implies Theorem~\ref{maintop}.
\vspace{2mm}

\noindent\textbf{Remark:} It is not 
necessary to work with chain complexes. The proof of Kontsetvich's
conjecture and all the major results of this paper have 
analogues in the category of spectra, which, however, 
are substantially more technically difficult.
\vspace{2mm}

\noindent\textbf{Acknowledgement:} I would like to thank Ralph 
Cohen, Igor Kriz and Sasha Voronov for important comments. 

\vspace{6mm}

\section{The Topological Story}
\label{topsec}

Let $X$ be a finite-dimensional simplicial complex, 
with a simplicial embedding in a sphere $S^N$ for some $N$ (for some
triangulation of $S^N$).  Instead of 
the normal bundle, we can construct
a parametrized spectrum $S(X)$ over $X$, similarly as the construction 
of the Spivak bundle for Poincare duality spaces.  
For each simplex $s$ of $X$, recall that $st(s)$, the star of $s$, is
defined to be the union of all open simplices $t$, such that $s$ is a 
face of $t$ (including $s$ itself). For each simplex $s$, 
consider the construction 
\begin{equation}
 (S^N \setminus (X \setminus st(s)))/(S^N \setminus X) . 
\label{quotient1}
\end{equation}
Here, by the symbol $/$, we mean the attachment of a cone on $S^N \setminus 
X$.  Note that the star is contravariantly functorial with respect to 
the simplicial structure: if $s_1$ and $s_2$ are simplices, and $s_1$ is a 
face of $s_2$, then $st(s_2) \subseteq st(s_1)$. Hence, the 
construction~(\ref{quotient1}) is also contravariantly
functorial with respect to the simplices of 
$X$. We define a simplicial space $S(X)_{\bullet}$ by 
\[ S(X)_n = \amalg_{dim(s) = n} (S^N \setminus (X \setminus (st(s)))/
(S^N \setminus X) \]
and define the parametrized spectrum $S(X)$ over $X$ to be the simplicial realization 
of $S(X)_{\bullet}$. To see that $S(X)$ is a parametrized spectrum
over $X$, note
that for each $n$-dimensional simplex $s$ of $X$, we have a map 
\begin{equation}
 ((S^N \setminus (X \setminus st(s)))/ (S^N \setminus X)) \times \Delta_n 
\rightarrow s 
\label{split1}
\end{equation}
by barycentric coordinates on $\Delta_n$. The sources of~(\ref{split1}) for 
all $n$-dimensional simplices $s$ make up $S(X)_n$, so taking disjoint unions 
and passing to the simplicial realization, we get a structure map 
\[ S(X) \rightarrow X . \]
Also, each map~(\ref{split1}) is split naturally with respect
to the simplicial structure, by going to the basepoint 
of $(S^N \setminus (X \setminus st(s)))/(S^N \setminus S)$ in the 
first coordinate, and taking the inverse of the canonical parametrization
$\Delta_n \rightarrow s$ in the second coordinate. So taking disjoint
unions and passing to the simplicial realization, we get a basepoint map 
\[ X \rightarrow S(X) . \] 

The parametrized spectrum $S(X)$ has the 
following property. Let $i: X \rightarrow \ast$ be the collapse map.
There is a functor $i_{\sharp}$ from spectra parametrized over $X$
to spectra, which is the left adjoint to the pullback functor $i^{\ast}$.
Specifically, $i_{\sharp}$ is obtained by collapsing the basepoint copy of 
$X$ to a single point~(see~\cite{hudual}). 
\begin{prop}
\[ i_{\sharp}S(X) \simeq \Sigma^N DX_{+} \]
where $DX_{+}$ is the
Spanier-Whitehead dual of $X_{+}$. 
\end{prop}
In particular, if $X$ is an orientable 
Poincare space, the $S(X)$ is homotopy equivalent
to the Spivak bundle of $X$. However,
for general $X$, $S(X)$ is not a spherical fibration over $X$.

\begin{proof}
For each simplex $s$ of $X$, we have an inclusion map 
\[ S^N \setminus (X \setminus st(s)) \rightarrow S^N . \]
Taking cones on $S^N \setminus X$ on both sides, and then applying 
simplicial realization, we get a map 
\begin{equation}
 i_{\sharp}S(X) \rightarrow S^N / (S^N \setminus X) . 
\label{quotient2}
\end{equation}
We will show that this is an equivalence. Note that simplicial realization 
and the cone are both colimits, so they commute with each other. Hence, it 
suffices to consider the simplicial realization 
$|S^N \setminus (X \setminus st(s))|$, and show that the map 
\[ |S^N \setminus (X \setminus st(s))| \rightarrow S^N \] 
is an equivalence. 
This is true since  
$(S^N \setminus (X \setminus st(s)))_{s}$ is 
an open cover of $S^N$, and the simplicial realization is its 
homotopy colimit. However, projection from the homotopy colimit of an open 
cover of a simplicial complex by open sets which are complements of 
subcomplexes to the original simplicial complex is an equivalence.
Hence,~(\ref{quotient2}) is an equivalence. But by 
Spanier-Whitehead duality, we also have that 
\[ S^N / (S^N \setminus X) \simeq \Sigma^N DX_+ . \]
\end{proof}

Consider the space $Maps(S^k, X)$ of all unbased continuous
maps $S^k \rightarrow X$. This is the higher dimensional 
analogue of the free loop space $LX$. Taking the evaluation map 
at a chosen basepoint of $S^k$
\[ ev: Maps(S^k, X) \rightarrow X \]
and pulling back $S(X)$, we can consider $S(X)$ as a parametrized 
spectrum over $Maps(S^k, X)$.  Also, there is the collapse map 
$j: Maps(S^k, X) \rightarrow \ast$. Pushing forward $S(X)$ along $j$
from paramectrized spectra over $Maps(S^k, X)$ to spectra (i.~e. applying 
the functor $j_{\sharp}$),
we get a spectrum $Maps(S^k, X)^{S(X)}$.  This is the generalization 
of the stable Thom complex $LM^{\nu_M}$.  

\vspace{5mm}

Now Deligne's conjecture, proven by McClure and Smith~\cite{ms}, states 
that the Hochschild cohomology complex of an associative algebra $R$
(i.~e. an algebra over $C_{\ast}\mathcal{C}_1$) has a natural structure of 
an algebra over $C_{\ast}\mathcal{C}_2$.  Kontsevich's conjecture calls for 
an analogue of this for algebras over $C_{\ast}\mathcal{C}_k$ for general $k$.
In~\cite{hkv}, we proved a version of Kontsevich's conjecture. Recall that 
if $R$ is an algebra over $C_{\ast}\mathcal{C}_k$, a $(C_{\ast}\mathcal{C}_k, R)$ 
module $M$
is given by structure maps 
\[ C_{\ast}\mathcal{C}_k(n+1) \otimes R^{\otimes n} \otimes M \rightarrow M \]
for each $n$, which satisfies the obvious associativity, equivariance, and 
unitality axioms. In~\cite{hkv}, we
defined closed model category structures on the 
category $C_{\ast}\mathcal{C}_k-Alg$ of $C_{\ast}
\mathcal{C}_k$-algebras, as well as 
on $(C_{\ast}\mathcal{C}_k, R)-Mod$, the category of modules over a 
$C_{\ast}\mathcal{C}_k$-algebra $R$. Namely, the fibrations and weak equivalences
are carried over from the category of chain complexes.
If $R$ is a $C_{\ast}\mathcal{C}_k$-algebra,
we defined the Hochschild cohomology complex of $R$ as 
\[ HC^{\ast}_{C_{\ast}\mathcal{C}_k}(R) = RHom_{(C_{\ast}\mathcal{C}_k, R)}(R, R) . \]

\begin{theor}[\cite{hkv}]
For a $C_{\ast}\mathcal{C}_k$-algebra $R$, $HC^{\ast}_{C_{\ast}\mathcal{C}_k}(R)$ 
has a natural structure as a $\mathcal{C}_{k+1}$-algebra.
\end{theor}

The analogue of Theorem~\ref{maintop} on 
chain complexes now follows from the following result.
\begin{theor}
For a $k+1$-connected simplicial complex 
$X$, there is a natural homotopy equivalence 
\[ C_{\ast}(Maps(S^k, X)^{S(X)}) \stackrel{\simeq}{\rightarrow} 
HC^{\ast}_{C_{\ast}\mathcal{C}_k}
(C^{\ast}X) . \]
\label{firsteq}
\end{theor}

\begin{proof}
The proof is analogous to that of~\cite{cj}. We consider $S^k_{\bullet}$, 
which is $S^k$ thought of as a simplicial
set. Then $Maps(S^k_{\bullet}, X)$ is a cosimplicial model for 
$Maps(S^k, X)$, whose $n$-th stage is given by 
$Maps(S^k_n, X)$. This induces in turn a cosimplicial decomposition of 
$Maps(S^k, X)^{S(X)}$. On the other hand, we have the following construction
due to Kelly~\cite{kelly}. Let $R$ be an $E_{\infty}$-algebra. We shall use the 
commutative, associative and unital product constructed in~\cite{km}. With 
respect to this product $\boxdot$, the $E_{\infty}$-structure map of $R$ is just
\[ R \boxdot R \rightarrow R \]
and it is strictly commutative, associative and unital. Then $\ast$ is the 
categorical coproduct of $E_{\infty}$-algebras in this language. Now 
following Kelly~\cite{kelly}, for a simplicial set $Z_{\bullet}$,
we can define the tensor product of $Z_{\bullet}$ and $R$ to 
be the simplicial realization
\[ Z_{\bullet} \otimes R = | \amalg_{Z_{\bullet}} R |  . \] 
Now note that the bar construction $B(R, R, R)$ can be thought of as 
$I_{\bullet} \otimes R$, where $I_{\bullet}$ is the unit interval with 
the standard simplicial 
structure. Analogously, for $k \geq 1$, consider a simplicial triangulation 
of $D^k$, with $S^k$ as a subcomplex. Let $M$ be a based $R$-module, 
(i.~e. we have a basepoint $K \rightarrow M$
to make $\boxdot$ defined, and a structure map 
\[ R \boxdot M \rightarrow M \]
with the usual commutative diagrams.) Then we have the higher 
bar construction 
\[ B^{(k)}(R, M) = (D^{k}_{\bullet} \otimes R) \boxdot_{S^{k-1}_{\bullet} 
\otimes R} M . \]
(Here, $M$ is an $(S^{k-1}_{\bullet} \otimes R)$-module via the collapse map 
$S^{k-1}_{\bullet} \otimes R \rightarrow R$.) We also have the higher 
cobar construction 
\[ C^{(k)}(R, M) = Hom_{S^{k-1}_{\bullet} \otimes R}(D^k_{\bullet} \otimes R, M) \]
where $Hom_{S^k_{\bullet} \otimes R}$ is the right adjoint to 
$\boxdot_{S^k_{\bullet} \otimes R}$. Now we have  

\begin{prop}
For an $E_{\infty}$-algebra $R$, 
\[ HC_{C_{\ast}\mathcal{C}_k}^{\ast}(R) \simeq C^{(k)}(R, R). \]
\label{highercob}
\end{prop}
 
Now to prove the theorem, note that evaluation gives a map 
\[ \Delta_n \times Maps(S^k_n, X)^{S(X)} \rightarrow 
(DX_{+})^{\wedge S^k_n} . \]
Passing to (appropriately rigidified) chains, and applying cosimplicial 
realization, we get a map 
\[ C_{\ast}(Maps(S^k, X)^{S(X)}) \rightarrow HC_{C_{\ast}\mathcal{C}_k}^{\ast}
(C^{\ast}X) . \]
An Eilenberg-Moore type spectral sequence then shows that this is an 
equivalence. 
\end{proof}

\begin{proof}[Proof of Proposition~\ref{highercob}]

For an $E_{\infty}$-algebra $R$, we have that 
$C^{(k)}(R, R)$ is defined to be 
$Hom_{S^{k-1}_{\bullet} \otimes R}(D^{k}_{\bullet} \otimes R, R)$. 
So it suffices to show that 
\[ S^{k-1}_{\bullet} \otimes R \stackrel{\simeq}{\rightarrow} A^{R}_1 \]
as algebras, such that the identity map $R \rightarrow R$ takes, up to 
quasi-isomorphism, the 
module structure on $R$ over $S^{k-1}_{\bullet} \otimes R$ to the 
module structure over 
$A_1^R$. We will show that this, in fact, follows from Poincare duality.

For an $k$-dimensional submanifold (possibly with boundary)
$M$ of ${\mathbb R}^k$ and a space $X$, we take 
\[ C^M_k (X) \]
to be the space of configurations of finitely many unordered 
little cubes in $M$, with faces orthogonal to the ${\mathbb R}^k$-coordinates, 
disjoint with the boundary of $M$, and 
labelled by elements of $X$. In particular, if $M = D^k$, then 
$C^{D^k}_k$ is quasi-isomorphic, as a right $C_k$-algebra, to $C_k$, 
the monad associated with the usual $k$-dimensional 
little cubes operad. For $M = S^{k-1} \times I$, the functor $C^{S^{k-1} \times I}_{k}$
is a left functor over $C_k$, and for an $\mathcal{C}_{\infty}$-algebra $R$
(considered as a $\mathcal{C}_k$-algebra by pullback), 
\begin{equation}
 A_1^R \simeq B(C_k^{S^{k-1} \times I}, C_k, R) . 
\label{barconst1}
\end{equation}
We claim that 
\begin{equation}
 B(C_k^{S^{k-1} \times I}, C_k , R) \simeq 
B(C_{\infty}^{S^{k-1} \times I \times I^{\infty}}, C_{\infty}, R) .
\label{barconst2}
\end{equation}
(Here, $(-)^{\infty}$ refers to $\colim_{n \rightarrow \infty} (-)^n$.)
To see this, one considers the double bar complex 
\begin{equation}
B((C_k^{S^{k-1} \times I}, C_k, C_{\infty}, C_{\infty}, R) .
\label{doublecx1}
\end{equation}
which is equivalent to $B(C_k^{S^{k-1} \times I}, C_k, R)$, and 
show that
\[ B(C_{k}^{S^{k-1} \times I}, C_k, C_{\infty}) \simeq 
C_{\infty}^{S^{k-1} \times I \times I^{\infty}} .   \]
Hence,~(\ref{doublecx1}) is 
also $B(C_{\infty}^{S^{k-1} \times I \times I^{\infty}}, C_{\infty}, R)$. 

Now it suffices to consider the case when $R = C_{\infty}X$ for a 
based space $X$. (Then 
\begin{equation*}
\begin{split}
B(C_{\infty}^{S^{k-1} \times I \times I^{\infty}}, C_{\infty}, R) & = 
B(C_{\infty}^{S^{k-1} \times I \times I^{\infty}}, C_{\infty}, C_{\infty}X) \\
& \simeq C_{\infty}^{S^{k-1} \times I \times I^{\infty}}X . )
\end{split}
\end{equation*}
We claim that this is equivalent to 
\[ (S^{k-1} \times I)_{\bullet} \otimes C_{\infty}X . \]

For $M$ as above, we also define $C^{\prime M}_k X$ to be 
the  space of unordered tuples of little cubes in $M$ with faces 
orthogonal to the ${\mathbb R}^k$-coordinates, labelled 
by elements of $X$, but differing from $C^{M}_k X$ in that 
little cubes are allowed to intersect with $\partial M$, but where we 
identify each configuration with the one obtained by deleting the 
little cubes which have nonempty intersections with $\partial M$. 
Let $M = \bigcup
U_i$ be a covering of $M$ by open neighborhoods $U_i$, where all intersections
are either empty or diffeomorphic to a convext set in ${\mathbb R}^k$. 
Then we have for $i_1, \ldots, i_n$ a natural equivalence
\[ C^{\prime U_{i_1} \cap \cdots \cap U_{i_n}}_k X \stackrel{\simeq}{\rightarrow}
Map(U_{i_1} \cap \cdots \cap U_{i_n}, \Sigma^k X) . \] 
Here, the suspension coordinate $\Sigma^k$ is equivalent to the restriction of 
the tangent bundle $\tau_M$ of $M$ to $U_{i_1} \cap \cdots \cap U_{i_n}$. 
Passing to inverse limits (in the \v{C}ech sense), we get that 
\begin{equation*}
\begin{split}
 C^{\prime M}_k X & \cong \ilim_{\leftarrow} 
C^{\prime U_{i_1} \cap \cdots \cap U_{i_n}}_k X \\
& \simeq \Gamma_M (S^{\tau_M} \wedge X) .
\end{split}
\end{equation*}
Here, $S^{\tau_M}$ is the sphere bundle of the tangent bundle of $M$, so 
$S^{\tau_M} \wedge X$ is a parametrized space over $M$, and $\Gamma_M$ 
denotes the space of global sections of this. Similarly, we also get that 
\[ C^M_k X \simeq \Gamma_{(M, \partial M)}(S^{\tau_M} \wedge X)   \]
where the right hand side denotes sections in which the boundary $\partial M$
goes to the basepoint. 
Applying this to $M \times I^r$ we get 
\begin{equation*}
\begin{split}
C_{k+r}^{M \times I^r} X & \simeq \Gamma_{(M \times D^r, M \times S^{r-1})} 
(S^{r-1} \wedge S^{\tau_M} \wedge X) \\
& \simeq \Gamma_M (\Omega^{r-1}(S^{r-1} \wedge S^{\tau_M} \wedge X)) .
\end{split}
\end{equation*}
Passing to $r \rightarrow \infty$, we get 
\begin{equation*}
\begin{split}
C_{\infty}^{M \times I^{\infty}}X & \simeq \Gamma_M (\Omega^{\infty}
(\Sigma^{\infty} S^{\tau_M} \wedge X)) \\
& \simeq \Omega^{\infty} Map(M^{\nu}, \Sigma^{\infty} X) \\
& \simeq \Omega^{\infty}\Sigma^{\infty}(M_{+} \wedge X). 
\end{split}
\end{equation*}
Here, $\nu$ denotes the stable normal bundle of $M$. The last equivalence 
is by Poincare duality. But this is also 
\[ \hocolim_{E_{\infty}} \Omega^{\infty}\Sigma^{\infty}((U_i)_{+} \wedge X) \]
where the homotopy colimit is taken in the category of 
$\mathcal{C}_{\infty}$-algebras. 
But we also have that
\begin{equation*}
\begin{split}
\Omega^{\infty}\Sigma^{\infty}((U_i)_{+} \wedge X) & \simeq 
U_i \otimes \Omega^{\infty} \Sigma^{\infty}X \\
& \simeq U_i \otimes C_{\infty}X . 
\end{split}
\end{equation*}
Here, the symbol $\otimes$ on the right hand side denotes simplicial 
realization of coproducts of $\mathcal{C}_{\infty}$-algebras in the 
sense of Kelly~\cite{kelly}.
The homotopy colimit is a \v{C}ech description of 
\[ M \otimes \Omega^{\infty}\Sigma^{\infty}X \simeq M \otimes C_{\infty}X . \]
Hence, we have an equivalence 
\begin{equation}
M \otimes C_{\infty}X \stackrel{\simeq}{\rightarrow} 
C_{\infty}^{M \times I^{\infty}} X . 
\end{equation}
This is, by definition, a map of right $C_{\infty}$-algebras 
(considering $X$ as a variable).
Now take $M$ to be a $(k-1)$-dimensional submanifold (without boundary)
of ${\mathbb R}^k$, and looking 
at $C^{M \times I \times I^{r}}_{k+r}X$, we get 
\begin{equation*}
\begin{split}
C_{k+r}^{M \times I \times I^r}X & \simeq \Gamma_{(M \times D^r, M \times S^{r-1}) 
\times (I, \delta I)} (S^{r-1} \wedge S^1 \wedge S^{\tau_M} \wedge X) \\
& \simeq \Gamma_M(\Omega^r(S^r \wedge S^{\tau_M} \wedge X) .
\end{split}
\end{equation*}
Taking $M = S^{k-1}$ and applying the above argument as $r \rightarrow 
\infty$, we get a map of right $\mathcal{C}_{\infty}$-algebras 
(with $X$ as a variable) 
\[ C_{\infty}^{S^{k-1} \times I \times I^{\infty}}X \simeq 
S^{k-1} \otimes C_{\infty}X . \]
Together with~(\ref{barconst2}) and~(\ref{barconst1}), we then have 
\begin{equation}
S^{k-1} \otimes R \stackrel{\simeq}{\rightarrow} B(C_{k}^{S^{k-1} \times I},
C_k, R) \simeq A_1^R 
\label{mainbarconst}
\end{equation}
as desired. (It is know that the 
$C_{\infty}$ Kelly $\otimes$-product corresponds
to the coproduct of structured $E_{\infty}$-algebras as described above; 
we omit the details here.)

To show that~(\ref{mainbarconst}) 
preserves the algebra structure, note that for 
$A^R_1$, the algebra structure is given by the gluing 
\[ (S^{k-1} \times I) \amalg (S^{k-1} \times I) \rightarrow 
S^{k-1} \times I . \]
Passing to the Poincare dual, this commutes with the $E_{\infty}$-structure coming 
from the $E_{\infty}$-structure of $R$, since $\mathcal{C}_1 \Box 
\mathcal{C}_{\infty} \simeq \mathcal{C}_{\infty}$, so an $A_{\infty}$-structure
in the category of $E_{\infty}$-algebras is absorbed into the 
$E_{\infty}$-structure. 

To show that the module structure on $R$ is preserved, note that 
$R = C_{\infty}X$ is a module over $A_1^R \simeq C_{\infty}^{S^{k-1} \times I 
\times I^{\infty}}X$ via the gluing 
\[ D^k \amalg (S^{k-1} \times I) \rightarrow D^k \]
using an equivalence 
\[ C_{\infty}X \simeq C_{\infty}^{D^k \times I^{\infty}} X \]
which is obtained similarly as above. This has the same naturality as the 
module structre of $R = C_{\infty}X$ over $S^{k-1}_{\bullet} \otimes 
C_{\infty}X$. 
\end{proof}


\vspace{3mm}

\section{The Algebraic Story}
\label{algsec}

In~\cite{hkv}, we outlined the following notion of Koszul dual derived categories,
using Quillen (co)homology~\cite{quillen}.
For a model category $C$ with products, let $Ab(C)$ be 
the full subcategory of abelian group objects (with respect to the 
categorical product) in $C$, and let $R: Ab(C) 
\rightarrow C$ be the forgetful functor.  (Given some 
set-theoretic conditions), this has a left adjoint
$L: C \rightarrow Ab(C)$, which is the abelianization functor. 
Let $F$ be the functorial cofibrant replacement in $C$, and let $L^{\prime} = 
LF$. Then $L^{\prime}$ has a right adjoint $R^{\prime}$. 
We define the (derived) 
Koszul transform $C_{!}$ of $C$ to be the category of coalgebras
over the comonad $L^{\prime}R^{\prime}$ on $Ab(C)$, and the Koszul dual category of 
$C$ to be $C^{!} = (C_!)^{op}$.  
For our purpose, what we will be interested in are the corresponding 
homotopy categories. To that end, we must define equivalences: equivalences 
in $Ab(C)$ are the maps $\varphi$ such that $R(\varphi)$ is an 
equivalence, and equivalences in $C_{!}$ are those morphisms which 
are equivalences in $Ab(C)$. As noted in~\cite{adams}, 
set-theoretical difficulties can arise in inverting equivalences in 
a given category. These difficulties can be overcome by introducing a closed 
model structure or another machinery. We do not know how to do this in 
the full generality, but in the particular cases considered below,
one easily sees that no set-theoretical difficulty arises
by small object arguments. 

We write 
$C^{!} = C^!_{Q}$ for this description of the Koszul dual category, to 
distinguish it from other descriptions of the Koszul dual. Here $Q$ stands for 
Quillen. Note that for an object $X$ of $C$,
$L^{\prime}X$ is in $(C_!)_Q$ via the unit of the adjunction 
$L^{\prime}X \rightarrow L^{\prime}R^{\prime}L^{\prime}X$.
We can therefore define
the Koszul dual functor $(-)^{!}_Q: C \rightarrow C^{!}_Q$ as 
the functor $L^{\prime}$.

Ginzburg and Kapranov~\cite{gk} studied Koszul duality for based operads
of chain complexes over a field $k$ of characteristic $0$ (based
just means that $\mathcal{C}(0) = K$).
Their definitions apply to that context only, but can be translated 
to say that the Koszul dual category of operads is itself, and 
for an operad $\mathcal{C}$, the Ginzburg-Kapranov 
Koszul dual operad $\mathcal{C}^!_{GK}$
is the free operad generated by the sequence of chain complexes 
given by 
\[ \mathcal{C}(n)^{\vee} = Hom(\mathcal{C}(n), K)   \]
with a suitable differential. In particular, an element of 
$\mathcal{C}(n)^{\vee}$ can be thought of as a directed tree, where
each vertex has at most one incoming edge, and there is exactly 
one vertex with no incoming edge (called the root); each vertex 
with one incoming edge and $n > 0$ outgoing edges is decorated by an 
element of $\mathcal{C}(n)^{\vee}$, the vertices with no outgoing edges
(called leaves) and the root are undecorated. Then the differential is obtained
by co-contraction of edges, i.~e. the operation dual to contracting 
a single edge in a tree. This is an analogue of the cobar 
construction for (co)-algebras. 

Note that the differential must have degree 1, so one gives each tree
a degree equal to its number of edges (the tree consisting of the root 
only has $-1$ edges, and the root is not counted as a leaf). 
For an $\mathcal{C}$-algebra $R$, one defines the Koszul 
dual $\mathcal{C}^{!}_{GK}$-algebra $R^{!}_{GK}$ as the free 
$\mathcal{C}^!_{GK}$-algebra on $\Sigma^{-1}R^{\vee}$, together with a
differential. Again, an element of $R^{!}_{GK}$ can be thought of 
as a directed tree with each vertex having at most one incoming edge, 
and exactly one root, where the inner vertices are decorated by elements of 
$\mathcal{C}(n)^{\vee}$, the root is undecorated, and the leaves are 
decorated by elements of $\Sigma^{-1}R^{\vee}$. The differential 
(of degree 1) is obtained by co-contracting a single edge where both 
vertices are decorated by elements of $\mathcal{C}^{\vee}$, or by 
co-contracting a ``bush'' where the root of the bush is decorated 
by an element of $\mathcal{C}(n)^{\vee}$ and the $n$ leaves are 
decorated by elements of $\Sigma^{-1}R^{\vee}$~(see~\cite{gk}). 

We remark that the reason for the shift in~\cite{gk} is that one 
wants to consider \emph{non-derived} Koszul duality, where one 
starts with an ungraded $\mathcal{C}$-algebra (where $\mathcal{C}$
is an operad in the category of $K$-modules), and obtain an ungraded 
non-derived Koszul dual algebra over the non-derived Koszul dual 
operad. This works for operads which~\cite{gk} call \emph{Koszul} 
(e.~g. the operad defining associative algebras is Koszul dual to 
itself; the operad defining commutative algebras is Koszul, and the non-derived
Koszul dual operad defines Lie algebras). Now in order for the 
ungraded non-derived Koszul duality to match the grading of the 
derived Koszul duality, the shift is needed. 

For an operad $\mathcal{C}$ over chain complexes, recall that 
we denoted by $\mathcal{C}-Alg$  
the category of algebras over $\mathcal{C}$, and that we defined 
a closed model structure on $\mathcal{C}-Alg$. We denote the
corresponding homotopy category by $h(\mathcal{C}-Alg)$.  

\begin{prop}
Let $\mathcal{C}$ be a based operad in the category of 
chain complexes of $K$-modules. Then the homotopy category of abelian 
objects in $\mathcal{C}$-algebras is equivalent to the homotopy category of 
chain complexes of 
$K$-modules. Also, there is an operad $\mathcal{C}^!_{cobar}$
such that the monad associated to $\mathcal{C}^!_{cobar}$ is equivalent to 
$(L^{\prime}R^{\prime})^{op}$. Moreover, the monad associated with the operad
$\mathcal{C}^!_{GK}$ is, via that equivalence and duality of 
chain complexes of $K$-modules, equivalent to the monad 
$(L^{\prime}R^{\prime})^{op}[1]$, obtained 
by first suspending an object, applying $L^{\prime}R^{\prime}$, then 
desuspending. 

Also, 
\begin{equation}
h(\mathcal{C}^!_{cobar}-Alg) \simeq h(\mathcal{C}-Alg)^{!}_{Q} .
\label{shiftbyone2}
\end{equation}
\label{gkkoszul} 
\end{prop}


\begin{proof}
To prove the first statement, note that if $X$ is an abelian 
$\mathcal{C}_k$-algebra, then $X \simeq K \oplus M$, where the 
$\mathcal{C}_k$-structure of $M$ is trivial. Therefore, the category 
of based abelian $C_k$-algebras is equivalent to the category of chain 
complexes, and this passes to the homotopy categories. 

We define $\mathcal{C}^!_{cobar}$ to be 
the free operad on $\{ \Sigma^{-1}\mathcal{C}(n)^{\vee} \}$, together 
with a differential. Namely, an element of $\mathcal{C}^!_{cobar}$ can be 
thought of as a directed tree, where each vertex has at 
most one incoming edge, and there is exactly one vertex with no 
incoming edge (the root); each vertex with 
one incoming edge and $n > 0$ 
outgoing edges is decorated by an element of $\Sigma^{-1} 
\mathcal{C}(n)^{\vee}$,
and the vertices with no outgoing edges (the leaves) and 
the root are undecorated. The differential 
has degree 1, and is obtained by co-contraction, i.~e. taking the dual of 
contracting a single edge of the tree. 

We can also consider $R^!_{cobar}$, the free $\mathcal{C}^!_{cobar}$-algebra
on $R^{\vee}$, together with a differential. An element of 
$R^!_{cobar}$ can be thought of as a directed tree where 
each vertex has at most one incoming edge, and where there is 
exact one root, where the vertices with 
one incoming edge and $n >0$ outgoing edges are decorated by elements of 
$\Sigma^{-1} \mathcal{C}(n)^{\vee}$, the root is 
undecorated, and the leaves are 
decorated by elements of $R^{\vee}$. The differential (of degree 1) is 
obtained by co-contracting a single edge where both vertices are decorated 
by elements of $\Sigma^{-1} \mathcal{C}^{\vee}$, or by co-contracting a 
``bush'', where the root of the 
bush is decorated by an element of $\Sigma^{-1} \mathcal{C}(n)^{\vee}$,
and the $n$ leaves are decorated by elements of $R^{\vee}$. In particular,
this is $\Sigma R^!_{GK}$, and we get that 
\begin{equation}
h(\Sigma (\mathcal{C}_{GK}^!-Alg)) \simeq h(\mathcal{C}_{cobar}^!-Alg) . 
\label{shiftbyone1}
\end{equation}
In fact, 
\[ \Sigma \mathcal{C}_{cobar}^! \simeq \mathcal{C}^!_{GK} . \]

To prove~(\ref{shiftbyone2}), for a 
$\mathcal{C}$-algebra $R$, we need to compare $R^!_{Q}$ and 
$R^!_{cobar}$. By a small object argument, 
we can assume that $R = CX$ for some $X$, where $C$ is 
the monad associated with $\mathcal{C}$. Then $R^!_Q \simeq HQ^{\ast}(R)$, 
the Quillen cohomology of $R$. This is known to be $X^{\vee}$. On the other 
hand, $R^!_{cobar}$ has the same shifts as in the cobar construction of 
monads $B(C, C, X)^{\vee}$ (by $-1$ on the vertices of the tree with outgoing 
edges, no shift on the copies of $R^{\vee}$ decorating
the leaves), so it is equivalent to 
$B(C, C, X)^{\vee} \simeq X^{\vee}$. Hence, $R^!_Q \simeq 
R^!_{cobar}$, and we get the equivalence~(\ref{shiftbyone2}). 
This uses the fact that an equivalence of monads induces an 
equivalence of homotopy categories of algebras where equivalences 
are induced by the forgetful functor; this can be proved by a small 
object argument.
 
\end{proof}

\noindent\textbf{Remark:}
As noted by the authors of~\cite{gk}, the reason Koszul duality for 
operads works is that the category of operads is Koszul dual to itself 
(with the machinery we set up, one can only state that for homotopy 
categories, but a non-derived statement is needed also for~\cite{gk}). 

Now this may be a puzzle, since operads involve symmetric group 
actions, and the operation $M/\Sigma_n$ is not self-dual. The explanation 
is that the derived category of abelian operads actually is not
equivalent to the derived category of sequences of chain complexes, 
but rather sequences $\{ M_n \}_{n \geq 1}$, where $M_n$ comes 
with a $\Sigma_n$-action. 
\vspace{3mm}

Now recall briefly Quillen cohomology of associative algebras. 
In the category of based (i.~e. augmented) associative algebras $R$ 
over a field $K$, abelian objects are just algebras of the form 
$K \oplus M$, where for $m_1, m_2 \in M$, $m_1m_2 = 0$. It 
follows that Quillen homology in this category is just the derived
functor of modules of indecomposables, which turns out to be 
the same as ``based Hochschild homology'' $Tor_{R}(K, K)$ with 
a shift of degrees by 1. Dualizing, we obtain a similar statement
for cohomology.


Kontsevich conjectured that analogously, Quillen and Hochschild 
cohomologies of based (i.~e. augmented) $C_{\ast}\mathcal{C}_k$-algebras
differ only by a shift of degrees by $k$. (By abuse of notation,
we shall write $\mathcal{C}_k$ instead of 
$C_{\ast}\mathcal{C}_k$ in the rest of the paper where there is no 
danger of confusion.) We prove this conjecture
here. By based Hochschild cohomology of a $\mathcal{C}_k$-algebra
$R$, we mean 
\[ HC^{\ast}_{based}(R) = RHom_{(\mathcal{C}_k, R)}(R, K) . \]
We grade the Hochschild and Quillen cohomology complexes homologically,
i.~e. with the higher cohomologies in negative degrees.

\begin{theor}
Let $R$ be a $\mathcal{C}_k$-algebra. Then for the Quillen and 
Hochschild cohomology complexes of $R$ with coefficients in $K$ satisfy
\[ HC^{\ast}_{based}(R) \simeq K \oplus HQ^{\ast}_{based}(R)[-k] . \]
(Here, the cohomology complexes are graded homologically, and $[-k]$
denotes shift down in dimension by $k$.)
\label{shift}
\end{theor}

We also recall briefly the suspension of operads. If $\mathcal{C}$ is 
an operad (in the category of chain complexes of $K$-modules), 
then the suspension $\Sigma \mathcal{C}$ is given by 
\[ (\Sigma \mathcal{C})(n) = \Sigma^{n-1} \mathcal{C}(n) . \]
One can check that this has the obvious operad structure. If $\mathcal{C}$
and $\mathcal{D}$ are operads, an equivalence of 
operads (with the usual model category structure where fibrations 
and equivalences are created by the forgetful functor to chain complexes)
\begin{equation}
  \Sigma^{-1} \mathcal{C} \simeq \mathcal{D}  
\label{suspequiv1}
\end{equation}
corresponds to a equivalence of homotopy categories 
\begin{equation}
 h(\mathcal{D}-Alg) \simeq h(\Sigma(\mathcal{C}-Alg)) . 
\label{suspequiv2}
\end{equation}
Here, $\Sigma(\mathcal{C}-Alg)$ denotes the category whose objects are of 
the form $\Sigma X$, where $X$ is a $\mathcal{C}$-algebra,
and whose morphisms are suspensions of morphisms of 
$\mathcal{C}$-algebras, and $h$ denotes the homotopy category. 
Namely, for any $\mathcal{C}$-algebra $X$, 
$\Sigma X$ has an obvious $\Sigma^{-1}\mathcal{C}$-action 
\[ \Sigma^{1-n}\mathcal{C}(n) \otimes 
(\Sigma X)^{\otimes n} \rightarrow \Sigma X . \]
This gives the equivalence~(\ref{suspequiv2}) given~(\ref{suspequiv1}).
Conversely, given~(\ref{suspequiv2}), comparing the
$\mathcal{C}$-action on $X$ 
\[ \mathcal{C}(n) \otimes X^{\otimes n} \rightarrow X \]
with the $\mathcal{D}$-action on $\Sigma X$
\[ \mathcal{D}(n) \otimes (\Sigma X)^{\otimes n} \rightarrow 
\Sigma X \]
we get 
\[ \Sigma^{n-1}\mathcal{D}(n) \simeq \mathcal{C}(n)   \]
for all $n$ by Proposition~\ref{monad2op} below. 

We also consider Koszul duality (in the Quillen sense)
for $\mathcal{C}_k$-algebras. 

\begin{theor}
There is an equivalence of derived ($=$ homotopy) categories
\[ h(\mathcal{C}_k-Alg)^!_Q \simeq h(\Sigma^{k} (\mathcal{C}_k-Alg)) . \]
\label{selfdual1}
\end{theor}

We also have the following result.
\begin{theor}
In the (derived) category of operads, 
\[ (\mathcal{C}_k)^{!}_{GK} \simeq \Sigma^{k-1}\mathcal{C}_k . \]
\label{selfdual2}
\end{theor}

We will use results of~\cite{hkv} in the proofs of Theorem~\ref{shift} and
Theorem~\ref{selfdual1}. However, \cite{hkv} was written in terms of 
simplicial categories, in particular the category of simplicial abelian 
groups. To reinterpret the results to the category of chain complexes, 
we will need the following lemma.

\begin{lemma}
There is a natural equivalence between the homotopy categories of 
$\mathcal{C}_k$-algebras in simplicial abelian groups and 
$\mathcal{C}_k$-algebras in connected chain complexes.
\label{catequiv}
\end{lemma}

Before proving the lemma, recall
that for an operad $\mathcal{C}$ in a closed symmetric monoidal category $Cat$,
there is a monad and a comonad associated
with it. Namely, let $\mathcal{C}-Alg$ be the subcategory of 
$\mathcal{C}$-algebras in $Cat$, and let $U: \mathcal{C}-Alg \rightarrow 
Cat$ be the forgetful functor. Then $U$ has a left adjoint 
$L: Cat \rightarrow \mathcal{C}-Alg$, and the monad $C$ associated 
with $\mathcal{C}$ is $RL$. Specifically, in the case when the 
category is that of topological spaces,
\[ CX = \amalg_{n \geq 0}(\mathcal{C}(n) \times_{\Sigma_n} X^n) . \]
A similar construction holds for chain complexes with the tensor
product, or for spectra with the smash product of spectra (more 
precisely $S$-modules~\cite{ekmm}).

Similarly, one can consider the category of coalgebras over 
an operad $\mathcal{C}$.
For an object $X$ of a closed symmetric monoidal category $Cat$, we can define 
the co-endomorphism operad $\overline{End}(X)$ of $X$ by 
\[ \overline{End}(X)(n) = \underline{Hom}(X, X^n) . \]
The operad structure maps are by symmetric monoidal products 
and compositions of internal $Hom$ objects.
Then $X$ is a coalgebra over $\mathcal{C}$ if there is a map of 
operads 
\[ \mathcal{C} \rightarrow \overline{End}(X) . \]
Equivalently, $X$ is a $\mathcal{C}$-coalgebra if there are maps 
\[ \mathcal{C}(n) \times X \rightarrow X^n  \]
or equivalently
\[ X \rightarrow \underline{Hom}(\mathcal{C}(n), X^n) \]
(in the language of spaces/simplicial sets) satisfying coassociativity, counitality, 
and equivariance axioms. Let $\mathcal{C}-Coalg$
denote the category of coalgebras over $\mathcal{C}$, and let 
$U: \mathcal{C}-Coalg \rightarrow Cat$ denote the forgetful functor. 
Then $U$ has a right adjoint $R$, and the comonad $\overline{C}$ associated 
with $\mathcal{C}$ is $UR$. Specifically, if $Cat$ is the category of 
spaces, we have 
\[ \overline{C}X = \prod_{n \geq 0}(\underline{Hom}(\mathcal{C}(n), 
X^n)^{\Sigma_n}) \]
which is the dual construction to $CX$.

\begin{proof}[Proof of Lemma~\ref{catequiv}]
The normalization functor from the category of simplicial abelian 
groups to the category of connected chain complexes is an equivalence
of categories, which does not preserve the tensor product structure. We
will compare the monad $C_k^{simp}$ of the operad $\mathcal{C}_k$ in 
simplicial abelian groups with the monad $C_k^{ch}$ of $\mathcal{C}_k$
in the category of connected chain complexes, and show that the equivalence
of categories induces an equivalence between the derived categories
of algebras over the 
monads in the two categories. 


First, the shuffle map gives a map of monads 
\[ C_k^{ch} \rightarrow C_k^{simp} \]
so any $C_k^{simp}$-algebra is a $C_k^{ch}$-algebra by pullback.
On the other hand, for a $C_k^{ch}$-algebra $R$, the 
2-sided bar construction of monads gives a simplicial $C_k^{simp}$-algebra
$B(C_k^{simp}, C_k^{ch}, R)$. But by Milnor's theorem, (the diagonal
simplicial realization of) a simplicial $C_k^{simp}$-algebra is 
a $C_k^{simp}$-algebra.
\end{proof}

\begin{proof}[Proof of Theorem~\ref{shift}]
We think of $\mathcal{C}_k$ as a based operad, and 
consider the category of based $\mathcal{C}_k$-algebras.
Let $R$ be a based $\mathcal{C}_k$-algebra.  
In~\cite{hkv}, we defined 
a monoid $A_1^R$, and showed that the homotopy category of 
$(\mathcal{C}_k, R)$-modules is equivalent to
the homotopy category of $A_1^R$-modules. The 
based Hochschild cohomology of $R$ with 
coefficients in $K$ is $RHom_{(\mathcal{C}_k, R)}(R, K)$, so it can be 
calculated as $RHom_{A_1^{R}-mod}(R, K)$. The based Hochschild homology of $R$ with 
coefficients in $K$ is its dual. In particular, it can be calculated
as the bar construction $B(K, A_1^R, R)$. 

On the other hand, there is also a suspension functor $\Sigma^k$ in 
the category of based $\mathcal{C}_k$-algebras. Namely, for a based 
$\mathcal{C}_k$-algebra $R$, $\Sigma^k R$ is defined to be the
the abelian $C_k$-algebra
\[ \Sigma^k R = (R \otimes C_{\ast} S^k) \oplus_{K \otimes C_{\ast} S^k} K . \]
(This functor is in fact defined for any chain complex, not just for 
a $\mathcal{C}_k$-algebra.)
Here $C_k$ is the monad associated with the operad $\mathcal{C}_k$. 
The functor $\Sigma^k$ is a right algebra over
the monad $C_k$, so we can also consider the bar construction of monads 
$B(\Sigma^k, C_k, R)$. We have the following lemma.

\begin{lemma}
There is a natural equivalence of right $C_k$-algebras
\[ B(K, A_1^{C_k}, C_k) \simeq \Sigma^k . \]
\label{suspension1}
\end{lemma}

Given Lemma~\ref{suspension1}, we consider the double bar complex 
\[ B(K, A_1^{C_k}, C_k, C_k, R) . \]
Looking at this as $B(B(K, A_1^R, C_k), C_k, R)$ and 
using Lemma~\ref{suspension1},
we see that it is $B(\Sigma^k, C_k, R)$. On the other hand, looking at it 
as $B(K, A_1^R, B(C_k, C_k, R))$, we get that it is also equivalent to 
$B(K, A_1^R, R)$. Hence, we get an equivalence
of $(\mathcal{C}_k, R)$-modules  
\begin{equation}
 B(\Sigma^k, C_k, R) \simeq B(K, A_1^R, R) 
\label{suspchange1}
\end{equation}
which is also the based Hochschild homology of $R$ with coefficients in 
$k$.

Now $B(\Sigma^k, C_k, R)$ is an abelian object in the category 
of $\mathcal{C}_k$-algebras. There is also a desuspension 
functor $\Omega^k$ from abelian $\mathcal{C}_k$-algebras to 
$\mathcal{C}_k$-algebras. 
Note that we have a natural map (not of algebras) in the derived category 
\begin{equation}
R \rightarrow \Omega^k \Sigma^k R
\label{loopsigma}
\end{equation}
(since $\Sigma^k R$ is equivalent to $K \oplus J[k]$, where $J$ is
the augmentation idea of $R$, and $[k]$ denotes shift by $k$).
Then~(\ref{loopsigma}) induces a map of 
$\mathcal{C}_k$-algebras 
\begin{equation}
C_k R \rightarrow \Omega^k \Sigma^k R .
\label{suspension2}
\end{equation}
In turn,~(\ref{suspension2}) induces a map of $\mathcal{C}_k$-algebras
\[ R \simeq B(C_k, C_k, R) \rightarrow B(\Omega^k\Sigma^k, C_k, R) .\]
The target is
the Hochschild homology of $R$ shifted by $k$, and is abelian. But now let 
$HQ_{\ast}(R)$ be the Quillen homology of $R$. As the derived functor of abelianization,
$HQ_{\ast}(-)$ is an universal functor from $\mathcal{C}_k$-algebras to abelian 
$\mathcal{C}_k$-algebras. Thus, we get a map of abelian $\mathcal{C}_k$-algebras
\begin{equation} 
HQ_{\ast}(R) \rightarrow B(\Omega^k \Sigma^k, C_k, R) . 
\label{suspchange2}
\end{equation}
This is an equivalence when $R = C_k X$ is a free $\mathcal{C}_k$-algebra. 
Hence, it is an equivalence for all $R$. Dualizing gives the equivalence 
between the Hochschild and Quillen cohomologies of $R$, up to a shift by $k$.

\end{proof}

\begin{proof}[Proof of Lemma~\ref{suspension1}]
We first consider the case when $X = C_{\ast}Y$ for a topological space 
$Y$. Then 
\[ B(K, A_1^{C_{k}X}, C_k X) \simeq C_{\ast} B(\ast, A_1^{C_k Y}, C_k Y) \]
and $\Sigma^k X \simeq C_{\ast} \Sigma^k Y$. So it suffices to show that 
\begin{equation}
B(\ast, A_1^{C_k Y}, C_k Y) \simeq \Sigma^k Y 
\end{equation}
on the level of spaces. 

Now $B(\ast, A_1^{C_k Y}, C_k Y)$ can be thought of 
(at least up to homotopy) as the homotopy colimit of 
the spaces of configurations 
of $n \geq 1$ nested concentric copies of $S^{k-1}$ in 
the interior of $D^k$, together with 
finitely many (unordered) points in $D^k$, which are decorated by elements of 
$Y$, and are not allowed to be on the copies of $S^{k-1}$, nor outside 
of the sphere with the largest radius. The homotopy colimit is formed by 
taking the nerve of the topological category whose object space is the 
disjoint union of the said configuration spaces, and morphisms are obtained
by omitting spheres. (When the sphere with the largest radius is omitted, 
any points falling outside the largest remaining sphere are deleted.)

On the other hand, let 
\begin{equation}
 B^{\prime}(\ast, A_1^{C_k Y}, C_k Y) 
\label{concentric1}
\end{equation}
be the homotopy colimit of configuration spaces defined in the same way,
with the only difference that this time we allow $Y$-decorated points in
$D^k$ outside of the sphere with the largest radius, but  a configuration 
is identified with the configuration obtained by deleting any points on 
$\partial D^k$. 
Then we have a forgetful map 
\begin{equation}
p: B^{\prime}(\ast, A_1^{C_k Y}, C_k Y) \rightarrow B(\ast, A_1^{C_k Y}, 
C_k Y)
\label{disappear1}
\end{equation}
by deleting the points outside the sphere with the largest radius.
Clearly, $p$ is a quasifibration and a homotopy equivalence. 


But now let the space $C_k^{\prime} Y$ be the space of configurations of 
finitely many (unordered) points in $D^k$, where a point on the boundary of $D^k$ 
disappears. It is know that 
\[ C_k^{\prime}Y \simeq Map(D^k, \Sigma^k Y) \simeq  \Sigma^k Y .  \]
We can define a map 
\begin{equation}
f: B^{\prime}(\ast, A_1^{C_k Y}, C_k Y) \rightarrow C^{\prime}_k Y
\label{disappear2}
\end{equation}
by forgetting all the concentric copies of $S^{k-1}$ inside $D^k$. We 
will show that~(\ref{disappear2}) is an equivalence. 

For each $r < 1$, let $U_r \subset C^{\prime}_k Y$ be the subspace
of configurations of points in $C_k^{\prime}Y$ such that there are no 
points on the copy of $S^{k-1}$ of radius $r$ (and the 
same center as $D^k$) inside $D^k$. The 
$\{ U_r \}_{r < 1}$ form an open covering of $C_k^{\prime}Y$, 
and we have  
\begin{equation}
 C_k^{\prime}Y \simeq \hocolim (U_{r_1} \cap \cdots \cap U_{r_n}) . 
\label{concentric2}
\end{equation}
where the homotopy colimit is formed in the usual \v{C}ech way (this
is a general fact). But now note that the right hand side 
of~(\ref{concentric2}) is homeomorphic to $B^{\prime}(\ast, 
A_1^{C_kY}, C_k Y)$. 



The above construction can be mimicked for 
general chain complex $X$ to produce a map of right 
$C_k$-algebras. More precisely, $B^{\prime}(\ast, 
A_1^{C_kX}, C_k X)$ can be defined by starting with~(\ref{concentric1})
for $Y = \ast$, forming the chain complex, tensoring each singular simplex with 
the appropriate number of copies of $X \oplus K$
(where $K$ is the base field), and identifying with
respect to permutations and deleting of points in the configurations. 
In this way, we obtain maps of right $C_{k}$-algebras
\begin{equation}
B(\ast, A_1^{C_k X}, C_kX) \simeq B^{\prime}(\ast, A_1^{C_k X}, C_k X) 
\rightarrow \Sigma^k X 
\label{concentric3}
\end{equation}
where $X$ is the variable. Now we know from the above discussion 
that these maps are equivalences when $X = C_k Y$. Since the 
constructions of~(\ref{concentric3}) clearly preserve equivalences, we 
know that this remains valid when $X$ is connected. To treat the 
general case, since we are dealing with modules over a field, we may 
assume the differential of $X$ is trivial. Then we may use the 
following trick: introduce an additional grading on $X$ where the second 
degree is always even. Then~(\ref{concentric3}) is bigraded. By choosing
the second degree suitably, we may assume that the total degree is 
$\geq 0$, and thus~(\ref{concentric3}) are equivalences. Therefore, the 
statement remains true when the second degree is forgotten. 
\end{proof}

We will need the following categorical 
lemma in the proof of Theorem~\ref{selfdual1}.

\begin{lemma}
Let $\overline{C}$ and $\overline{D}$ be comonads, and suppose 
that $\overline{C}$ coacts on every $\overline{D}$-coalgebra, such that 
the coaction is a natural transformation in the category 
of $\overline{D}$-coalgebras. Then 
there is a natural map of comonads $\overline{\alpha}: \overline{D} \rightarrow 
\overline{C}$.
\label{comonad1}
\end{lemma}

\begin{proof}
It suffices to consider the statement for monads: let $C$ and $D$ be 
monads, and suppose that $C$ acts naturally on every $D$-algebra. Then 
there is a natural map of monads $\alpha: C \rightarrow D$. Let $X$ be any object.
We define the map $\alpha$ on $X$ to be the composition
\[ \alpha_{X} : CX \stackrel{C\eta_D}{\rightarrow} CDX \stackrel{\gamma}{\rightarrow}
DX  \]
where $\eta_D$ is the unit of the monad $D$, and $\gamma$ is the action of 
$C$ on $D$. So $\alpha$ is a natural transformation. We need to show that $\alpha$ is 
a map of monads. The unitality axiom follows from the following diagram, 
which commutes by the naturality of $\eta_C$: 
\[ \diagram
X \rto^{\eta_D} \dto_{\eta_C} & DX \dto_{\eta_C} \drto^{=} & \\
CX \rto_{C\eta_D} & CDX \rto_{\gamma} & DX .
\enddiagram \]
We also need to consider the diagram 
\[ \diagram 
CCX \rto^{CC\eta_D} \dto_{\mu_C} & CCDX \drto_{\mu_C} \rto^{C\gamma} & 
CDX \rto^{C\eta_D} & CDDX \rto^{\gamma} & DDX \dto^{\mu_D} \\
CX \rrto_{C\eta_D} && CDX \rrto_{\gamma} && DX .
\enddiagram \]
Here, $\mu_C$ and $\mu_D$ are the multiplication structures of the 
monads $C$ and $D$ respectively.
By naturality, the left and right squares of this diagram commute, so the 
large rectangle also commutes. 
\end{proof}

\begin{prop}
Let $\mathcal{C}, \mathcal{D}$ be operads in the category of 
chain complexes of $K$-modules, and let $C$, $D$ be the associated
monads. Suppose that $C$, $D$ are equivalent, i.~e. there exists a diagram 
of monads of the form 
\[ \diagram 
& \dlto_{\simeq} E \drto^{\simeq} & \\
C & & D 
\enddiagram \]
where $\simeq$ denotes equivalence. Then $\mathcal{C}$, $\mathcal{D}$ 
are equivalent operads. 
\label{monad2op}
\end{prop}

\begin{proof}
Let $nK = K \oplus \cdots \oplus K$. Then $\mathcal{C}(n)$ is a 
wedge summand of $C(nK)$, which is the cokernel (not cofiber) of the map 
\begin{equation}
\oplus_{\varphi_i} C((n-1)K) \rightarrow C(nK) 
\label{coker1}
\end{equation}
where the direct sum is over all injections $\varphi_i: \{ 1, \ldots, n-1 \} 
\rightarrow \{ 1, \ldots, n \}$. One can check that for any monad $C$,
the cokernel of~(\ref{coker1}) is a direct summand and naturally an operad. 
One way to do this is as follows: we shall define, by induction, a canonical
summand $C^{\prime}(nK)$ of $C(nK)$ such that 
\begin{equation}
C(nK) = \oplus_{S \subseteq \{1, \ldots, n\}} C^{\prime}(|S|K). 
\label{summand1}
\end{equation}
where~(\ref{summand1}) is given by the retraction
\begin{equation}
C^{\prime}(|S|K) \subseteq C(|S|K) \subseteq C(nK) \rightarrow 
C(|S|K) \rightarrow C^{\prime}(|S|K)
\label{summand2}
\end{equation}
for each $S$. (Here, the two middle maps are induced by the inclusion 
$S \subseteq \{ 1, \ldots, n \}$ and the projection $nK \rightarrow |S|K$
which sends $\{ i \}K \rightarrow \{ i \}K$ by the identity for $i \in S$,
and $\{ i \}K \rightarrow 0$ for $i \not\in S$.) 

Suppose~(\ref{summand1}) holds with $n$ replaced by $n-1$. Then consider the maps
\begin{equation}
\oplus_{S \subset\neq \{ 1, \ldots, n\}}C^{\prime}(|S|K) \rightarrow 
C(nK) \rightarrow \oplus_{S \subset\neq \{ 1, \ldots, n \}} C^{\prime}
(|S|K) 
\label{summand3}
\end{equation}
defined by~(\ref{summand2}). Then the composition~(\ref{summand3}) is 
the identity by the induction hypothesis. Let $C^{\prime}(nK)$ be the 
canonical direct complement of the retraction~(\ref{summand3}). The 
induction hypothesis follows.

Then $C^{\prime}$ is the cokernel of~(\ref{coker1}), and the monad properties
of $C$ translate under the splitting~(\ref{summand1}) to saying 
that $C(nK)$ is an operad. Thus,
the statement follows from the fact that a direct summand of an equivalence
is an equivalence. 
\end{proof}

\begin{proof}[Proof of Theorem~\ref{selfdual1}]
Recall that in the proof of Theorem~\ref{shift}, we 
established (see~(\ref{suspchange1})) that the based Quillen cohomology 
$HQ^{\ast}_{based}(R)$ is naturally equivalent to the bar 
construction $B(\Sigma^k, C_k, R)$. 

However, the operad
$\mathcal{C}_k$ coacts on the $k$-dimensional sphere.
Concretely, given an element $x$ of $\mathcal{C}_k(n)$, we have a map 
\[ f_x: S^k \rightarrow (S^k)^n \]
as follows. The element $x$ is a configuration of $n$ $k$-dimensional 
disks inside a single $k$-dimensional disk. Let 
$f_x =(f_1, \ldots, f_n)$. The map $f_i: S^k \rightarrow S^k$
is defined by letting
the $i$-th little disk in the domain map to the target by scaling, 
and letting the rest of the domain go to the basepoint. 
The maps $f_x$ specify a map
\[ \mathcal{C}_k(n) \rightarrow Hom(S^k, (S^k)^n) \]
which is a coaction of $\mathcal{C}_k$ on $S^k$.  This gives 
rise to a coaction of the comonad $\overline{C}_k$ associated 
with $\mathcal{C}_k$ on the functor $\Sigma^k$, so 
$\mathcal{C}_k$ coacts on the Quillen cohomology $B(\Sigma^k, C_k, R)$.


Denote by $L^{\prime}$ the functor $B(\Omega^k\Sigma^k, C_k, -)$ from based 
$\mathcal{C}_k$-algebras to abelian based $\mathcal{C}_k$-algebras
(recall that the category of abelian based $\mathcal{C}_k$-algebras
is equivalent to the category of chain complexes over $K$). This 
functor has a right adjoint $R^{\prime}$. Then $R^{\prime}$ is equivalent to 
the forgetful functor by~(\ref{suspchange2}). Moreover, (\ref{suspchange2})
shows that the comonad $L^{\prime}R^{\prime}$ is equivalent to 
$L \cdot F\cdot R$ where $R$
is the forgetful functor from abelian based $\mathcal{C}_k$-algebras
to based $\mathcal{C}_k$-algebras, $L$ is its left adjoint, and $F$
is cofibrant replacement. 
Therefore, the Koszul transform of the derived category of 
based $\mathcal{C}_k$-algebras is equivalent to the derived category of 
$L^{\prime}R^{\prime}$-algebras. 
But the arguments made earlier in this proof show that
the comonad $\overline{C_k}$ associated with $\mathcal{C}_k$ coacts 
on the $k$-suspension of
every $L^{\prime}R^{\prime}$-algebra.
By Lemma~\ref{comonad1}, we obtain a map of 
comonads
\begin{equation}
\overline{\alpha}: L^{\prime}R^{\prime} \rightarrow \overline{C_k}[-k] 
\label{comonad2}
\end{equation}
where $\overline{C_k}[-k]$ means the comonad obtained by
suspending an object by $k$, applying $\overline{C_k}$, and 
desuspending by $k$. 

It suffices to prove that the map $\overline{\alpha}$ is an equivalence. 
Now by construction, the 
map $\overline{\alpha}$ is 
\begin{equation}
B(S^k, C_k, X) \rightarrow \overline{C_k} B(S^k, C_k, X) 
\rightarrow \overline{C_k} \Sigma^k X .
\label{comonad3}
\end{equation}
We would like to have control over the cohomology of~(\ref{comonad3}). 
We have the following lemma. 

\begin{lemma}
Let $Y$ be an $1$-connected space. Then the 
$\mathcal{C}_k$-algebra $C^{\ast}\Sigma^k Y$ 
is equivalent to an abelian $\mathcal{C}_k$-algebra. 
\label{connected1}
\end{lemma}

Given Lemma~\ref{connected1}, we assume that $X = C^{\ast}\Sigma^k Y$ 
for some $1$-connected space $Y$. Then by passing to the Poincare-dual 
description of Proposition~\ref{highercob}, we get that
\begin{equation}
 H^{\ast}(B(S^k, C_k, X)) =H^{\ast}\Omega^k \Sigma^k Y 
\label{control1}
\end{equation}
and on cohomology, the first map of~(\ref{comonad3}) is the dual to the 
Dyer-Lashof operations. The point is that iterating the bar 
construction on an abelian $\mathcal{C}_k$-algebra 
$K \oplus M$ gives objects whose cohomologies are free 
graded-commutative algebras on Steenrod operations, which, 
if $M$ is in negative homological degrees, correspond 
precisely to Dyer-Lashof operations. We also have that 
\[ H^{\ast}(\Sigma^k X) \simeq H^{\ast}(Y) \]
and the second map of~(\ref{comonad3}) in cohomology is projection to 
the characteristic class
\[ H^{\ast}\Omega^k \Sigma^k Y \rightarrow H^{\ast} Y \]
induced by the unit of the adjunction 
\[ Y \rightarrow \Omega^k \Sigma^k Y . \]
Passing to the dual formulations, one gets that the 
composition~(\ref{comonad3}), using the 
identification~(\ref{control1}), induces the identity on cohomology. So it is 
an quasi-isomorphism.
Now we can show that~(\ref{comonad3}) is an equivalence for a 
general $X$ by the same double grading trick as applied at the end 
of the proof of Lemma~\ref{suspension1}.

\end{proof}

\begin{proof}[Proof of Lemma~\ref{connected1}]
We first consider the case when $k=1$. The
$\mathcal{C}_1$-algebra (i.~e. $A_{\infty}$-algebra) structure on 
$C^{\ast}S^1$ is induced by the diagonal map 
\[ \Delta : S^1 \rightarrow S^1 \wedge \cdots \wedge S^1 . \]
Let $(S^1, \Delta)$ denote $S^1$ with the $A_{\infty}$-coalgebra structure 
given by $\Delta$, and let $(S^1, 0)$ denote $S^1$ with the trivial 
$A_{\infty}$-coalgebra structure. Recall also that the loop of the 
space of $A_{\infty}$-coalgebra structures on $X$ is 
$End_{A_{\infty}}(X, X)$ (with the chosen $A_{\infty}$-structure on $X$ 
as the basepoint). We would like to construct a map 
\begin{equation}
(S^1, \Delta) \rightarrow (S^1, 0) . 
\label{trivialstr1}
\end{equation}
For this, we have a spectral sequence 
\begin{equation}
RHom_{A_{\infty}-alg}(\Lambda[x], \Lambda[x]) \Rightarrow 
RHom_{A_{\infty}-alg}((C^{\ast}S^1, \Delta), (C^{\ast}S^1, 0)) 
\label{trivialstr2}
\end{equation}
where $x$ is a generator in (topological) dimension 1, and the second
$\Lambda[x]$ is abelian as an $A_{\infty}$-algebra (i.~e. has trivial
multiplication). But we also have that 
\[ RHom_{A_{\infty}-alg}(\Lambda[x], \Lambda[x]) \cong 
HQ^{\ast}_{Asso}(\Lambda[x], K) . \]
By the grading of the spectral sequence, the differentials of~(\ref{trivialstr2})
have no nontrivial targets, so the element sof 
$HQ^{\ast}_{Asso}(\Lambda[x], K)$ survive and we get that for $S^1$, 
$C^{\ast}S^1$ is an abelian $A_{\infty}$-algebra. 

For general $k$, the $\mathcal{C}_k$-algebra structure on 
$C^{\ast}S^k$ is induced by the diagonal 
\[ \Delta : S^k \rightarrow S^k \wedge \cdots \wedge S^k . \]
But this is the smash product of $k$ copies of 
\[ \Delta : S^1 \rightarrow S^1 \wedge \cdots \wedge S^1 \]
and the $\mathcal{C}_k$-coaction on $S^k$ breaks up as the $\Box$ 
of $k$ copies of the coaction of $\mathcal{C}_1$ on $S^1$, since 
$\mathcal{C}_k \simeq \mathcal{C}_1 \Box \cdots \Box \mathcal{C}_1$. 
Hence, the $k=1$ case gives that the $\mathcal{C}_k$-action on 
$C^{\ast}S^k$ is equivalent to the trivial one. 
\end{proof}

Finally, we prove Theorem~\ref{selfdual2}.

\begin{proof}[Proof of Theorem~\ref{selfdual2}]

By~(\ref{comonad2}) and Proposition~\ref{gkkoszul}, the 
corresponding monads are equivalent. Now use Proposition~\ref{monad2op}.
\end{proof}

\vspace{3mm}



\begin{proof}[Proof of Proposition~\ref{cochain}]
We use Theorems~\ref{shift} and~\ref{selfdual1} and their proofs. 
The Koszul dual of an object is its Quillen cohomology, which in this case
is the shift of the dual of its based Hochschild homology. Therefore, 
the Koszul dual of $C_{\ast}\Omega^k X$ is the dual of 
\[ B(\Sigma^k, C_k, C_{\ast}\Omega^k X) \simeq C_{\ast}
B(\Sigma^k, C_k, \Omega^k X) \simeq C_{\ast}X . \]
This is $C^{\ast}X$. (Recall the shift at the end of the proof of 
Theorem~\ref{selfdual1}; that proof also gives control over the 
$\mathcal{C}_k$-algebra structure, and one can show that the 
$\mathcal{C}_k$-algebra structure of $C^{\ast}X$ one gets coincides
with the $E_{\infty}$-structure one gets by the diagonal. Namely, by 
naturality, the $\mathcal{C}_k$-algebra structure commutes with the 
$E_{\infty}$-structure; therefore, both structures 
coincide by Lemma~\ref{boxeinfty} below. )
\end{proof}
\vspace{6mm}

\section{A Duality for Hochschild Cohomology}
\label{hochsec}

Given Theorem~\ref{selfdual1}, we can compare the Hochschild cohomology
complex
of a $C_{\ast}\mathcal{C}_k$ algebra $R$ with that of its Koszul dual 
$R^!$.  By~\cite{hkv}, both are algebras over 
$C_{\ast}\mathcal{C}_{k+1}$.

\begin{theor}
For a connected (with respect to the homological grading)
$C_{\ast}\mathcal{C}_k$-algebra $R$, such that $H_{\ast}(R^!)$ 
is a finite-dimensional
$K$-module, there is a natural equivalence 
of chain complexes
\[ HC^{\ast}_{C_{\ast}\mathcal{C}_k}(R) \simeq HC^{\ast}_{C_{\ast}\mathcal{C}_k}
(R^!) . \]
\label{kossame}
\end{theor}

\noindent\textbf{Remark:}
One may conjecture that~(\ref{kossame}) is true as $\mathcal{C}_{k+1}$-algebras,
but I do not have have a proof at this point.

\begin{lemma}
Let $F^nC$, $n \geq 0$ be a complete decreasing filtration of 
a chain complex $C$, so that $gr_{F} C$ is acyclic. Then $C$ is 
acyclic. (By ``complete'', we mean that $\cap_n F^n C =0$, and
if $(x_n)$ is a sequence of elements in 
$C$, such that $x_n - x_{n-1} \in F^n C$, $n \rightarrow \infty$, then 
there exists an element $x$ such that $x_n -x \in F^m C$ for some $m$, and 
$m \rightarrow \infty$.)
\label{acyclic}
\end{lemma}

\begin{proof}
Suppose that $dx = 0$. Then this is true in $C/F^1C$. Then there 
exists an element $y_0 \in C$, such that
\[ x_1 = dy_0 -x \in F^1 C  . \]
Hence, $dx_1 =0$, so $dx_1 = 0 \in F^1C/F^2C$. So there exists
$y_1 \in F^1C$, such that 
\[ x_2 = dy_1 - x_1 \in F^2 C. \]
Analogously, we construct $x_n$, $y_n$ such that 
\[ x_n = dy_{n-1} -x_{n-1} \in F^n C . \]
Let $y = \sum y_i$. Then $dy = x$. 
\end{proof}

We shall also need the following lemma, which is proved analogously. (Note
that cosimplicial realization is a limit, so the filtration by degrees 
is complete.)

\begin{lemma}
Let $\varphi: C_{\bullet \bullet} \rightarrow D_{\bullet \bullet}$ be a 
map of double cosimplicial complexes, which induces an equivalence 
\[ \varphi_n: |C_{\bullet}|_n \rightarrow |D_{\bullet}|_n \]
for each $n \geq 0$. Then 
\[ \varphi: |C_{\bullet \bullet}| \rightarrow |D_{\bullet \bullet}| \]
is an equivalence.
\label{3dim}
\end{lemma}

We shall first give a complete proof of 
Theorem~\ref{kossame} in the case $k =1$. This case is 
in some sense classical and the proof in this case is more direct,
and it motivates the methods used in the general case. A reader 
not interested in this case can skip directly to the case of 
general $k$. 

\begin{proof}[Proof of Theorem~\ref{kossame} for the case $k=1$:]
In the case of $k = 1$, we can assume $R$ is an associative 
algebra. A model of 
the Koszul transform $R_{!}$ in this case is given by the bar construction
$B(R)$. Thus, 
in this case, $R^!$ can also be described as the dual of 
the bar construction $B(K, R, K)^{\vee}$, where $K$ is the base field.
Here, as above, $(-)^{\vee}$ denotes $Hom_K(-, K)$. 
We will compute both $HC^{\ast}(R)$ and $HC^{\ast}(R^!)$
by the following construction. 
We construct the double cosimplicial chain complex
whose $(m,n)$-th term is
\begin{equation}
 Hom_{K}(R^{\otimes m}, B(R, R, K) \otimes B(K, R, K)^{\otimes n} \otimes 
B(K, R, R)) .
\label{double1}
\end{equation}
The $m$-variable cosimplicial structure coincides with 
\begin{equation}
Hom_{R \otimes R^{op}}(B(R, R, R), M_n)
\label{homcosimp1}
\end{equation}
where $M_n$ is the $(R \otimes R^{op})$-module
\[ B(R, R, K) \otimes B(K, R, K)^{\otimes n} \otimes B(K, R, R) . \]
On the other hand, the $n$-cosimplicial structure is induced from 
the cosimplicial structure of 
\[ M_{\bullet} = C(B(R, R, K), B(K, R, K), B(K, R, R)) . \]
It is easy to check that both cosimplicial structures commute. 

First, we shall compare~(\ref{double1}) with $HC^{\ast}(R^{!})$. To this 
end, note that $HC^{\ast}(R^!)$ is equivalent to the cobar 
complex 
\begin{equation}
Hom_{(B(R) \otimes B(R)^{op})-comodules}(B(R), C(B(R), B(R), B(R)))
\label{barkos1}
\end{equation}
by the hypothesis that $R^! \simeq B(R)^{\vee}$ is homologically finite.
Now recalling~(\ref{homcosimp1}), the cosimplicial
$(B(R) \otimes B(R)^{op})$-comodule 
\[ Hom_{R \otimes R^{op}}(B(R, R, R), M_0) \]
is equivalent to the $(B(R) \otimes B(R)^{op})$-comodule
\[ Hom_{(B(R) \otimes B(R)^{op})-comodules}(B(R), B(R) \otimes B(R)^{op}) . \]
(Both are models of the dual of $B(R)$ in the category of 
comodules over the coalgebra $B(R) \otimes B(R)^{op}$.)
But then by Lemma~\ref{3dim}, the bicosimplicial realization of~(\ref{double1})
and~(\ref{barkos1}) are equivalent, i.~e. to
\[ HC^{\ast}(R^!) . \]

To show that~(\ref{double1}) is equivalent to $HC^{\ast}(R)$, by 
Lemma~\ref{3dim}, it suffices to show that we have an equivalence
of $(R \otimes R^{op})$-modules
\begin{equation}
\begin{split}
& C(B(R, R, K), B(K, R, K), B(K, R, R)) \\
& = |B(R, R, K) \otimes 
B(K, R, K)^{\otimes \bullet} \otimes B(K, R, R)| \\
& \simeq R.
\end{split}
\label{barkos2}
\end{equation}
But consider the canonical map 
\begin{equation}
B(R, R, K) \Box_{B(K, R, K)} B(K, R, R) \rightarrow 
C(B(R, R, K), B(K, R, K), B(K, R, R)).
\label{barkos3}
\end{equation}
Then the left hand side is isomorphic to $B(R, R, R)$, so it suffices to 
prove that~(\ref{barkos3}) is an equivalence. But consider 
the double decreasing (bounded below) filtration of~(\ref{barkos3})
obtained by filtering the first and last copies of $R$ on both 
sides by dimension. (The subset $R_{\geq n} \subset R$ is a submodule, 
so $B(R_{\geq n}, R, K)$ is a sub-$B(R)$-comodule of $B(R, R, K)$.)
The filtration is complete, because in a fixed total degree
a sequence of elements of increasing dimensions must have increasing
cohomological degrees. Then the associated graded map is of the form 
\[ M \Box_{B(R)} N \rightarrow C(M, B(R), N) \]
where $M, N$ are extended $B(R)$-comodules, which is always an 
equivalence. Since the filtration is complete and bounded below, the 
map~(\ref{barkos3}) is an equivalence by Lemma~\ref{acyclic}.
\end{proof}

\begin{proof}[Proof of Theorem~\ref{kossame} for general $k$:]
We will make 
use of the monoid $A_1 = A_1^{R}$ constructed in~\cite{hkv}. 
For a $\mathcal{C}_k$-algebra $R$, the derived category of $(\mathcal{C}_k,
R)$-modules is equivalent to the derived category of $A_1^R$-modules. 
$A_1^R$ is constructed as the free $(\mathcal{C}_k, R)$-module on 
one generator. Now the 
analogue of~(\ref{double1}) will be the double cosimplicial object whose
$(m, n)$-th term is
\begin{equation}
Hom_{K}((A_1^R)^{\otimes m} \otimes R, B(A_1^R, A_1^R, K) \otimes 
B(K, A_1^R, K)^{\otimes n} \otimes B(K, A_1^R, R)) .
\label{doublek}
\end{equation}
The $m$-cosimplicial structure is described as the 
cosimplicial structure of 
\[ Hom_{A_1^R}(B(A_1^R, A_1^R, R), M_n) \]
where 
\[ M_n = B(A_1^R, A_1^R, K) \otimes B(K, A_1^R, K)^{\otimes n}
\otimes B(K, A_1^R, R) . \]
The $n$-cosimplicial structure of~(\ref{doublek}) is described as 
\[ Hom_{K}((A_1^R)^{\otimes m} \otimes R, C(B(A_1^R, A_1^R, K), 
B(K, A_1^R, K), B(K, A_1^R, R))). \]
By arguments similar to that 
of the case when $k = 1$, we have that 
\[ C(B(A_1^R, A_1^R, K), B(K, A_1^R, K), B(K, A_1^R, R)) 
\simeq B(A_1^R, A_1^R, R) \simeq R \]
as $A_1^R$-modules. Hence,~(\ref{doublek}) is equivalent 
to 
\[ Hom_{A_1^R}(B(A_1^R, A_1^R, R), R) \]
which is $HC^{\ast}_{\mathcal{C}_k}(R)$, the Hochschild cohomology complex of 
$R$ in the category of $\mathcal{C}_k$-algebras, analogously to the 
second part of the proof for $k =1$. 
By the proof of 
Theorem~\ref{shift}, we have that 
$B(K, A_1^R, R) \simeq B(\Sigma^k, C_k, R)$, which 
is a description for the $k$-fold internal suspension of $R$ in the 
category of $\mathcal{C}_k$-algebras. This is the Koszul transform 
$R_!$ of $R$ in $\mathcal{C}_k$-algebras, so we can describe the 
Koszul dual $R^!$ as $B(K, A_1^R, R)^{\vee}$. 

To show that~(\ref{doublek}) also gives the Hochschild cohomology complex 
of $R^!$, we will show that
\begin{equation} 
B(K, A_1^R, K) \simeq (A_1^{R^!})^{\vee}
\label{dualize}
\end{equation}
as coalgebras over $\mathcal{C}_k$, and that the
canonical coaction of $(A_1^{R^!})^{\vee}$ on $R_!$ 
under~(\ref{dualize}) corresponds to the canonical coaction of
$B(K, A_1^R, K)$ on $B(K, A_1^R, R)$. Then, analogously to the first
part of the proof of the $k= 1$ case,~(\ref{doublek}) is equivalent to
\begin{equation}
\begin{split}
Hom_{B(K, A_1^R, K)-comodules} ( & B(K, A_1^R, R), \\
& C(B(K, A_1^R, K), B(K, A_1^R, K), B(K, A_1^R, R)))
\end{split}
\label{comodhom1}
\end{equation}
because
\[ Hom_{B(K, A_1^R, K)-comodules}(B(K, A_1^R, R), B(K, A_1^R, K)) \]
and 
\[ Hom_{A_1^R}(B(A_1^R, A_1^R, R), B(A_1^R, A_1^R, K)) \]
are equivalent (right) $B(K, A_1^R, K)$-comodules. Now~(\ref{comodhom1}) is 
a model of $HC^{\ast}_{\mathcal{C}_k}(R^!)$ by~(\ref{dualize}) (convergence
discussions are the same as in the $k =1$ case.)

To prove~(\ref{dualize}) (with its addendum about the 
coaction on $B(K, A_1^R, R)$), 
it suffices to show it when $R$ is a free $\mathcal{C}_k$-algebra, i.~e.
when $R = C_k(Y)$ for some chain complex $Y$.  In this 
case, by Lemma~\ref{suspension1}, we have that 
\[ B(K, A_1^{C_k(Y)}, C_k(Y)) \simeq B(\Sigma^k, C_k, C_k(Y)) 
\simeq \Sigma^k Y  . \]
To control the $B(K, A_1^R, K)$-comodule structure, by the 
double-grading trick in the proof of Lemma~\ref{suspension1},
we can assume that 
$Y \simeq C_{\ast}(X)$ for some $(k+1)$-connected topological space $X$. 
Now dually to the notion of a 
$(\mathcal{C}_k, R)$-module for a $\mathcal{C}_k$-algebra $R$, for 
a coalgebra $R^{\prime}$ over $\mathcal{C}_k$, we have the notion of 
a $(\mathcal{C}_k, R^{\prime})$-comodule $N$, whose structure maps are of the 
form 
\[ \mathcal{C}_k(n+1) \otimes N \rightarrow (R^{\prime})^{\otimes n} 
\otimes N  \]
with the obvious axioms.
Applying $(-)^{\vee}$ takes a $\mathcal{C}_k$-algebra to 
a $\mathcal{C}_k$-coalgebra, and vice versa. It also takes a 
$(\mathcal{C}_k, R)$-module to a comodule over $(\mathcal{C}_k, 
R^{\vee})$, and vice versa. 
Dual to the free module description 
of $A_1^R$, we have that given a coalgebra $S$ over $\mathcal{C}_k$, 
$(A_1^{S^{\vee}})^{\vee}$ is the cofree $(\mathcal{C}_k, S)$-comodule on one generator. 
So the right hand side of~(\ref{dualize}) can be described as 
the cofree $(\mathcal{C}_k, \Sigma^k Y)$-comodule on one generator
(because $\mathcal{C}_k^{!} \simeq \mathcal{C}_k$). We will
show that the left hand side of~(\ref{dualize}) is a configuration space 
model of the same. 
We have that $B(K, A_1^{C_kY}, K)$ is equivalent to the chain complex
of a semi-simplicial (without degeneracies) topological space $Z$, 
whose $n$-th stage is the space of all configurations of $n+1$ 
distinct $(k-1)$-spheres $s_0, \ldots, s_n$ in ${\mathbb R}^k$ with 
center $0$, where the radius of $s_i$ is less the radius of $s_{i+1}$, and
the $k$-dimensional annulus between $s_i$ and $s_{i+1}$ contains an 
undetermined number of (unordered) disjoint $k$-disks decorated by elements of $X$.

However, now $Z$ is equivalent to the configuration space of 
unordered tuples of finitely many disjoint points in $S^{k-1} \times I$
which are decorated by elements of $X$, and are not allowed to 
collide, but a point disappears when on the boundary. But this is in turn 
a configuration space model for 
\begin{equation}
Maps(S^{k-1} \times I, \Sigma^k X) .
\label{configsp1}
\end{equation}

But now let $\mathcal{D}(n)$ be the space of all ordered $n$-tuples
of disjoint $k$-dimensional disks in $S^{k-1} \times I$, and for a space
$X$ let 
\begin{equation*}
\begin{split}
 \underline{D}_0(X) = & \{ (f_n)_n \in \prod_n Maps(\mathcal{D}(n), 
X^n)^{\Sigma_n} \ |\ \\
& \mathrm{for\ the\ degeneracies\ } \epsilon_i: \mathcal{D}(n) \rightarrow 
\mathcal{D}(n-1), e_i: X^{n-1} \rightarrow X^n, f_n \epsilon_i = e_i f_n   \} .
\end{split}
\end{equation*}
Similarly, define
\begin{equation*}
\begin{split}
 & \underline{C_k}_{, 0}(X) = \{ (f_n)_n \in \prod Maps(\mathcal{C}_k(n), 
X^n)^{\Sigma_n} \ |\  \\
& \mathrm{for\ the\ degeneracies\ } \epsilon_i: \mathcal{D}(n) \rightarrow 
\mathcal{D}(n-1), e_i: X^{n-1} \rightarrow X^n, f_n \epsilon_i = e_i f_n  \} .
\end{split}
\end{equation*}
Then one can show that~(\ref{configsp1}) is equivalent to the 
cobar construction of spaces 
\begin{equation}
Cobar(\underline{D}_0, \underline{C_k}_{, 0}, \Sigma^k X) .
\label{configsp2})
\end{equation}
An Eilenberg-Moore type spectral sequence can further be used
to commute $C_{\ast}$ past the cobar construction, to show that 
the chain complex of~(\ref{configsp2}) is equivalent to the 
cofree $(\mathcal{C}_k, \Sigma^k C_{\ast} X)$-comodule on one generator.
\end{proof}




Finally, we have the following lemma. 

\begin{lemma}
For every operad $\mathcal{D}$ where $\mathcal{D}(1) \simeq \ast$, 
$\mathcal{D} \Box \mathcal{C}_{\infty}$ is an $E_{\infty}$-operad. 
\label{boxeinfty}
\end{lemma}

\begin{proof}
We need to verify that each $(\mathcal{D} \Box \mathcal{C}_{\infty})(n)$ 
is contractible. We will first prove a weaker statement, constructing 
a homotopy between the inclusion 
\begin{equation}
\mathcal{D} \rightarrow \mathcal{D} \Box \mathcal{C}_{\infty}
\end{equation}
and a constant map. 
Choose $n \geq 0$ and $\gamma \in \mathcal{C}_{\infty}(n)$, and 
let $h_i$ be paths between $1$ and $\gamma(\ast, \ldots, 1, \ldots, \ast)$
where $1$ is in the $i$-th place, $i = 1, \ldots, n$. Now let 
$x \in \mathcal{D}(n)$. Consider the homotopy which on $x$ acts as the 
path 
\begin{equation}
\begin{split}
x = x (1, \ldots, 1) & \simeq x (\gamma(1, \ast, \ldots, \ast), \ldots, 
\gamma(\ast, \ldots, \ast, 1)) \\
& = \gamma (x(1, \ast, \ldots, \ast), \ldots, x(\ast, \ldots, \ast, 1)) .
\end{split}
\end{equation}
Here, the first equivalence is via $x (h_1, \ldots, h_n)$. Since 
$\mathcal{D}(1)$ is contractible, the right hand side of this is homotopic 
to $\gamma(1, \ldots, 1) = \gamma$. 

Now consider any map 
\begin{equation}
\alpha: K \rightarrow \mathcal{D} \Box \mathcal{C}_{\infty}
\label{einftymap1}
\end{equation}
where $K$ is compact. Then by the definition of 
$\mathcal{C}_{\infty}$ (with the weak topology), (\ref{einftymap1}) 
factors through 
\[ K \stackrel{\alpha_k}{\rightarrow} \mathcal{D} \Box \mathcal{C}_{k} 
\stackrel{\iota}{\rightarrow} \mathcal{D} \Box \mathcal{C}_{\infty} \]
for some $k$. Here, the second map $\iota$ is the inclusion map. However, 
note that then by~\cite{hkv}, we have a commutative diagram 
\[ \diagram 
\mathcal{D} \Box \mathcal{C}_k \rrto^{\iota} \drto_{\subseteq} && 
\mathcal{D} \Box \mathcal{C}_{\infty} \\
 & \mathcal{D} \Box \mathcal{C}_{k} \Box \mathcal{C}_{\infty} 
\urto_{Id_{\mathcal{D}} \Box \pi} & 
\enddiagram \]
where $\pi$ is the composition map $\mathcal{C}_{k} \Box 
\mathcal{C}_{\infty} \rightarrow \mathcal{C}_{\infty}$. By the previous 
argument (applied to $\mathcal{D} \Box \mathcal{C}_{k}$ in place of 
$\mathcal{D}$), the map $\subseteq$ of the diagram is null-homotopic. Hence, 
so is $\subseteq \cdot \alpha_k$, and so is $\iota \cdot \alpha_k = 
Id_{\mathcal{D}} \Box (\subseteq \cdot \alpha_k)$. 
\end{proof}

\end{document}